\documentclass{article}

\usepackage{amssymb, amsmath, latexsym, mathrsfs, graphicx, stmaryrd,
amsthm, psfrag, color, multirow}
\usepackage{rotating}
\usepackage{url}
\usepackage[enableskew, vcentermath]{youngtab}
\newtheorem{thm}{Theorem}[section]
\newtheorem{lemma}[thm]{Lemma}
\newtheorem{cor}[thm]{Corollary}

\newtheorem{conjecture}[thm]{Conjecture}

\theoremstyle{definition}

\newtheorem{definition}[thm]{Definition}

\newcommand{\Z}{\mathbb{Z}}

\newcommand{\R}{\mathbb{R}}
\newcommand{\C}{\mathbb{C}}

\makeatletter
\def\revddots{\mathinner{\mkern1mu\raise\p@
\vbox{\kern7\p@\hbox{.}}\mkern2mu
\raise4\p@\hbox{.}\mkern2mu\raise7\p@\hbox{.}\mkern1mu}}
\makeatother
\def \abs#1{\lvert #1 \rvert}

\DeclareMathOperator{\sgn}{sgn}
\renewcommand{\epsilon}{\varepsilon}

\DeclareMathOperator{\cone}{cone}

\begin{document}

\title{On the Computation of Clebsch--Gordan Coefficients and the
       Dilation Effect
       \footnote{Research supported by NSF Grants DMS-0309694 and
                 VIGRE Grant DMS-0135345.  The first author is also
                 grateful for support received from the Alexander von
                 Humboldt foundation and a UC Davis Chancellor's
                 fellow award.
       }
}

\author{ Jes\'us A. De Loera and Tyrrell B. McAllister }

\maketitle

\vskip .5cm

\begin{abstract}
    \noindent We investigate the problem of computing tensor product
    multiplicities for complex semisimple Lie algebras.  Even though
    computing these numbers is \( \#P \)-hard in general, we show that
    if the rank of the Lie algebra is assumed fixed, then there is a
    polynomial time algorithm, based on counting the lattice points in
    polytopes.  In fact, for Lie algebras of type $A_r$, there is an
    algorithm, based on the ellipsoid algorithm, to decide when the
    coefficients are nonzero in polynomial time for arbitrary rank.
    Our experiments show that the lattice point algorithm is superior
    in practice to the standard techniques for computing
    multiplicities when the weights have large entries but small rank.
    Using an implementation of this algorithm, we provide experimental
    evidence for conjectured generalizations of the saturation
    property of Littlewood--Richardson coefficients.  One of these
    conjectures seems to be valid for types \( B_n \), \( C_n \), and
    \( D_n \).
\end{abstract}

\section{Introduction}

Given highest weights \( \lambda \), \( \mu \), and \( \nu \) for a
finite dimensional complex semisimple Lie algebra, we denote by \(
C^\nu_{\lambda \mu} \) the multiplicity of the irreducible
representation \( V_{\nu} \) in the tensor product of \( V_{\lambda}
\) and \( V_{\mu} \); that is, we write
\begin{equation}\label{eq:tensormult}
    V_{\lambda} \otimes V_{\mu} 
    =
    \bigoplus_{\nu} C_{\lambda \mu}^{\nu} V_{\nu}.
\end{equation}
In general, the numbers \( C_{\lambda \mu}^{\nu} \) are called
\emph{Clebsch--Gordan coefficients}.  In the specific case of type \(
A_{r} \) Lie algebras, the values \( C_{\lambda \mu}^{\nu} \) defined
in equation (\ref{eq:tensormult}) are called {\em
Littlewood--Richardson coefficients}.  When we are specifically
discussing the type \( A_{r} \) case, we will adhere to convention and
write \( c_{\lambda \mu}^{\nu} \) for \( C_{\lambda \mu}^{\nu} \).

The concrete computation of Clebsch--Gordan coefficients (sometimes
known as the {\em Clebsch--Gordan problem} \cite{fultonharris}) has
attracted a lot of attention from not only representation theorists,
but also from physicists, who employ them in the study of quantum
mechanics (e.g. \cite{cohengraaf, wybourne}).  The importance of these
coefficients is also evidenced by their widespread appearance in other
fields of mathematics besides representation theory.  For example, the
Littlewood--Richardson coefficients appear in combinatorics via
symmetric functions and in enumerative algebraic geometry via Schubert
varieties and Grassmannians (see for instance \cite{stanley3,fulton}).
More recently, Clebsch--Gordan coefficients are playing an important
role on the study of $P$ vs.  $NP$ (see \cite{mulmuleysohoni}).  Very
recently, Narayanan has proved that the computations of
Clebsch--Gordan coefficients is in general $\#P$-complete
\cite{narayanan}.  Here are our contributions:

\noindent {\bf (1)} We combine the lattice point enumeration algorithm
of Barvinok \cite{B:94} with polyhedral results due to Knutson and Tao
\cite{knutsontao} and Berenstein and Zelevinsky \cite{BZ:01} in the
polyhedral realization of Clebsch--Gordan coefficients to produce a
new algorithm for computing these coefficients. Our main
theoretical result is:

\begin{thm}\label{main}
    For fixed rank \( r \), if $\mathfrak{g}$ is a complex semisimple
    Lie algebra of rank \( r \), then one can compute Clebsch--Gordan
    coefficients for $\mathfrak{g}$ in time polynomial in the input
    size of the defining weights.
    
    Moreover, as a consequence of the polynomiality of linear
    programming and the saturation property of Lie Algebras of type
    $A_r$, deciding whether \( c_{\lambda \mu}^{\nu} = 0 \) can be
    done in polynomial time, even when the rank is not fixed.
\end{thm}

\noindent{\bf (2)} We implemented the algorithm
for types \( A_{r} \), \( B_{r} \), \( C_{r} \), and \( D_{r} \) (the
so-called ``classical'' Lie Algebras) using the software packages
\texttt{LattE} and \texttt{Maple}.  In many instances, our implementation
performs faster than standard methods, such as those implemented in
the software \texttt{LiE}. Our software is freely available at
\url{http://math.ucdavis.edu/~tmcal}.

\noindent {\bf (3)} Via computer experiments, we explored general properties
satisfied by the Clebsch--Gordan coefficients for the classical Lie
algebras under the operation of {\em stretching of multiplicities} in
the sense of \cite{kingetal}.  On the basis of abundant experimental
evidence, we propose two conjectured generalizations of the Saturation
Theorem of Knutson and Tao \cite{knutsontao}.  One of them, which
applies to all of the classical root systems, is an extension of
earlier work by King et al.  (\cite{kingetal}).

\noindent Organization of the paper: In Section \ref{computation1}, after a
review of the background material, we prove Theorem \ref{main}.
Section \ref{comparison} explains our experiments comparing our
software, a mixture of \texttt{Maple} and \texttt{LattE} \cite{latte},
with \texttt{LiE}.  In Section \ref{conjectures}, we present the two
conjectures, both of which, if true, would generalize the Saturation
Theorem of Knutson and Tao.

\section{Clebsch--Gordan coefficients: Polyhedral Algorithms}
\label{computation1}

As stated in the introduction, we are interested in the problem of
efficiently computing $C_{\lambda \mu}^{\nu}$ in the tensor product
expansion \( V_{\lambda} \otimes V_{\mu} = \bigoplus_{\nu} C_{\lambda
\mu}^{\nu} V_{\nu} \).  It appears that the most common method used to
compute the Clebsch--Gordan coefficients is based on Klimyk's formula (see
lemma below).  For example, it is used in \texttt{LiE} \cite{lie} and
the \texttt{Maple} \cite{maple} packages \texttt{Coxeter} and
\texttt{Weyl} \cite{coxeterweyl}.  

\begin{lemma} (\cite{H:72}, Exercise. 24.9)
    Fix a complex semisimple Lie algebra \( \mathfrak{g} \), and let
    \( \mathfrak{W} \) be the associated Weyl group.  For each weight
    \( \nu \) of \( \mathfrak{g} \), let \( \sgn(\nu) \) denote the
    parity of the minimum length of an element \( \sigma \in
    \mathfrak{W} \) such that \( \sigma(\nu) \) is a highest weight,
    and let \( \{\nu\} \) denote that highest weight.  Let \( \delta
    \) be one-half the sum of the positive simple roots of \(
    \mathfrak{g} \).  Finally, for each highest weight \( \lambda \)
    of \( \mathfrak{g} \), let \( K_{\lambda \nu} \) be the
    multiplicity of \( \nu \) in \( V_{\lambda} \).
    
    Then, given highest weights \( \lambda \) and \( \mu \) of \(
    \mathfrak{g} \), we have that
    \[ 
        V_{\lambda} \otimes V_{\mu}
        =
        \bigoplus_{\epsilon}
            K_{\lambda \epsilon} \sgn(\epsilon + \mu + \delta)
            V_{\{\epsilon + \mu + \delta \} - \delta},
    \]
    where the sum is over weights \( \epsilon \) of \( \mathfrak{g} \) 
    with trivial stabilizer subgroup in \( \mathfrak{W} \).
\end{lemma}

Implementations of Klimyk's algorithm begin by computing the weight
spaces appearing with nonzero multiplicity in the representation
$V_\lambda$.  Then, for each such weight $\epsilon$ with trivial
stabilizer, one computes the Weyl group orbit of \( \epsilon + \mu +
\delta \).  One then finds the dominant member of the orbit and notes
the number \( l \) of reflections needed to reach it.  Finally, one
adds \( (-1)^l K_{\lambda \epsilon} \) to the multiplicity of \(
V_{\{\epsilon + \mu + \delta \} - \delta} \).

From the point of view of computational complexity, Klimyk's algorithm
has two main disadvantages.  First, it requires the computation of
weight space multiplicities, which is in general a \( \#P \)-hard
problem \cite{narayanan}.  The second disadvantage is that the
algorithm requires visiting all of the orbit members, which can be an
exponentially large set.  Indeed, Klimyk's formula above is
exponential in the size of the input weights. Thus, in practice, the
sizes of $\lambda$, $\mu$, and $\nu$ usually need to be small.  One
can then ask for an algorithm that behaves well as the sizes of the
input weights increase, at least if some other parameter is fixed.
Stembridge also raised the challenge of crafting algorithms based on
geometric ideas such as Littelmann's paths or Kashiwara's crystal
bases \cite{kashiwara,littelmann} (see comment on page 29, section 7,
of \cite{stembridge}).  As we see below, there is such an algorithm,
based on the polyhedral geometry of the Clebsch--Gordan coefficients.

In 1992, Berenstein and Zelevinsky presented a combinatorial
interpretation of the Littlewood--Richardson coefficients as the
number of lattice points in members of a certain family of polytopes
\cite{berensteinzelevinskyI}.  In 1998, Knutson and Tao introduced
another family, the \emph{hive polytopes}, which they used to prove
the Saturation Theorem (see \cite{knutsontao} and Section
\ref{conjectures}).  Each of the polytopes presented by Berenstein and
Zelevinsky in 1992 is the image under an injective lattice-preserving
linear map of a hive polytope \cite{pakvallejo}.  Therefore \(
c^\nu_{\lambda \mu} \) equals the number of integer lattice points in
a corresponding hive polytope \( H_{\lambda \mu}^\nu \).  Finally, in
2001, Berenstein and Zelevinsky \cite{BZ:01} introduced polytopes
which enumerate Clebsch--Gordan coefficients for any finite
dimensional complex semisimple Lie algebra.  We refer to this last
family of polytopes as the {\em BZ--polytopes}.  We now give the exact
definition of the hive polytopes.  These polytopes exist in the
polyhedral cone of \emph{hive patterns}, which we now define.

\begin{definition}
    Fix \( r \in \Z_{\geq 0} \) and let \( \mathcal{H} = \{ (i,j,k)
    \in \Z^{3}_{\geq 0} : i + j + k = r \} \).  A {\em hive pattern}
    is a map
    \[
        h \colon \mathcal{H} \rightarrow \R_{\geq 0},
        \qquad 
        (i,j,k) \mapsto h_{ijk},
    \]
    satisfying the {\em rhombus inequalities}:
    \begin{itemize}
        \item \( h_{i,j-1,k+1} + h_{i-1,j+1,k} \leq h_{ijk} + h_{i-1,j,k+1}, \)
        
        \item \( h_{ijk} + h_{i-1,j-1,k+2}     \leq h_{i,j-1,k+1} + h_{i-1,j,k+1}, \)
        
        \item \( h_{i+1,j-1,k} + h_{i-1,j,k+1} \leq h_{ijk} + h_{i,j-1,k+1}. \)
    \end{itemize}
    for \( (i,j,k) \in \mathcal{H} \), \( i, j \geq 1 \).
\end{definition}

Equivalently, a hive pattern of rank \( r \) is a triangular array of
real numbers
\[
    \begin{matrix}
                  &           &             &             & h_{0,0,r}   &             &             &        &           \\
                  &           &             &             &             &             &             &        &           \\
                  &           &             & h_{1,0,r-1} &             & h_{0,1,r-1} &             &        &           \\
                  &           &             &             &             &             &             &        &           \\
                  &           & h_{2,0,r-2} &             & h_{1,1,r-2} &             & h_{0,2,r-2} &        &           \\
                  &           &             &             &             &             &             &        &           \\
                  & \revddots &             &             & \vdots      &             &             & \ddots &           \\
                  &           &             &             &             &             &             &        &           \\
        h_{r,0,0} &           & h_{r-1,1,0} &             & \cdots      &             & h_{1,r-1,0} &        & h_{0,r,0}
    \end{matrix}
\]
such that, in every ``little rhombus'' of entries
\[ 
    \begin{matrix}
    c &   & b &   \\
      &   &   &   \\
      & a &   & d
        \end{matrix} \qquad, \qquad
    \begin{matrix}
      & c &   \\
      &   &   \\
    a &   & b \\
      &   &   \\
      & d &
    \end{matrix} \qquad, \qquad
    \begin{matrix}
      & a &   & c \\
      &   &   &   \\
    d &   & b &   
    \end{matrix}
\]
we have \( a+b \geq c+d \).
                                                                       
Here is example of a hive pattern with \( r = 5 \):
\[
    \begin{array}{ccccccccccc}
           &    &    &    &    & 0  &    &    &    &    &   \\
           &    &    &    & 5  &    & 8  &    &    &    &   \\
           &    &    & 8  &    & 12 &    & 13 &    &    &   \\
           &    & 11 &    & 15 &    & 17 &    & 18 &    &   \\
           & 12 &    & 16 &    & 18 &    & 20 &    & 20 &
    \end{array}
\]

Recall that when \( \mathfrak{g} \) is of type \( A_{r} \), so that \(
\mathfrak{g} \cong \mathfrak{sl}_{r+1}(\C) \) for some \( r \geq 2 \),
the highest weights are (with respect to the canonical basis)
partitions of length \( r \), \textit{i.e.}, sequences \( \lambda \)
of integers \( \lambda_{1} \geq \dotsb \geq \lambda_{n} \geq 0 \).  We
write \( \abs{\lambda} \) for \( \sum_{i} \lambda_{i} \), the {\em
size} of the partition \( \lambda \).

\begin{definition} \label{def:HivePolytopes}
    Given partitions \( \lambda, \mu, \nu \in \Z^{n}_{\geq 0} \), the
    {\em hive polytope} \( H_{\lambda \mu}^{\nu} \) is the set of hive
    patterns with boundary entries fixed as below.
    \[
        \begin{array}{ccccccccccccc}
                                                         &           &         &         &         & 0       &         &         &         &        &         \\
                                                & & & & \makebox[0pt][r]{\( \nu_{1} = \) } \bullet &         & \bullet \makebox[0pt][l]{ \( = \lambda_{1} \)} & & & & \\
                              & & & \makebox[0pt][r]{\( \nu_{1} + \nu_{2} = \) } \bullet &         & \bullet &         & \bullet \makebox[0pt][l]{ \( = \lambda_{1} + \lambda_{2} \)} & & & & \\
            & & \makebox[0pt][r]{\( \nu_{1} + \nu_{2} + \nu_{3} = \) } \bullet &         & \bullet &         & \bullet &         & \bullet \makebox[0pt][l]{ \( = \lambda_{1} + \lambda_{2} + \lambda_{3}\)} & & & \\
                                                         & \revddots &         &         &         & \vdots  &         &         &         & \ddots &         \\
            \makebox[0pt][r]{\( \abs{\nu} = \) } \bullet &           & \cdots  &         & 
                                                              \rotatebox[origin=c]{-45}{\( \bullet \) \makebox[0pt][l]{\( = \abs{\lambda} + \mu_{1} + \mu_{2} + \mu_{3} \)}}
                                                                                                   &         &
                                                                                  \rotatebox[origin=c]{-45}{\( \bullet \) \makebox[0pt][l]{\( = \abs{\lambda} + \mu_{1} + \mu_{2} \)}}
                                                                                                                       &         &
                                                                                                      \rotatebox[origin=c]{-45}{\( \bullet \) \makebox[0pt][l]{\( = \abs{\lambda} + \mu_{1} \)}} 
                                                                                                                                           &         & \bullet \makebox[0pt][l]{ \( = \abs{\lambda} \)}
        \end{array}
    \]
\end{definition}

\bigskip
\bigskip
\bigskip
\bigskip
\bigskip
\bigskip

Note that, for fixed \( r \), the input size of a hive polytope \(
H_{\lambda \mu}^{\nu} \) grows linearly with the input sizes of the
weights \( \lambda, \mu \), and \( \nu \).  We will need the following
result:

\begin{lemma}\label{hives} 
    (\cite{KTW:?}) The Littlewood--Richardson coefficient \(
    c_{\lambda \mu}^{\nu} \) equals the number of integer lattice
    points in the hive polytope \( H_{\lambda\mu}^{\nu} \).
\end{lemma}

Unfortunately, the description of the BZ-polytopes is more involved
than that of the hive polytopes above.  Therefore, we refer the reader
to Theorems 2.5 and 2.6 of \cite{BZ:01}, which give their description
as systems of linear equalities and inequalities in terms of the root
systems \( B_{r} \), \( C_{r},\) and \( D_{r} \).  The reader may also
view the contents of our maple notebook, available at
\url{http://math.ucdavis.edu/~tmcal}, for completely explicit
descriptions of the necessary inequalities.  The specific properties
of the BZ-polytopes that we need to prove our theorem are (1) for
fixed rank \( r \), the dimensions of the BZ-polytopes are bounded
above by a constant, (2) the input size of a BZ-polytope grows
linearly with the input sizes of the weights \( \lambda \), \( \mu \),
and \( \nu \), and (3) the following result describing the
relationship between BZ-polytopes and Clebsch--Gordan coefficients:

\begin{lemma} \label{bz}
    (Theorem 2.4 of \cite{BZ:01}) Fix a finite dimensional complex
    semi-simple Lie algebra \( \mathfrak{g} \) and a triple of highest
    weights $(\lambda, \mu, \nu)$ for \( \mathfrak{g} \).  Then the
    Clebsch--Gordan coefficient \( C_{\lambda \mu}^{\nu} \) equals the
    number of integer lattice points in the corresponding BZ-polytope.
\end{lemma}

Finally, the last ingredient necessary is A. Barvinok's algorithm for
counting lattice points in polytopes in polynomial time for fixed
dimension.  Several detailed descriptions of the algorithm in Lemma
\ref{barvithm} are now available the literature (see
\cite{deloeraetal} and all the references therein).

\begin{lemma} \label{barvithm}
    (\cite{B:94}) Fix \( d \in \Z_{\geq 0} \).  Then, given a system
    of equalities and inequalities defining a rational convex polytope
    \( P \subset \R^{d} \), we can compute \( \#(P \cap Z^{d}) \) in
    time polynomial in the input size of the polytope.
\end{lemma}

\noindent {\em Proof of Theorem \ref{main}:}
    First, if we fix the rank of the Lie algebra, then we fix an upper
    bound on the dimension of the hive or BZ polytope.  Moreover, the
    input sizes of these polytopes grow linearly with the input sizes
    of the weights.  Thus, by Barvinok's theorem (Lemma \ref{barvithm}
    stated above), their lattice points can be computed in time
    polynomial in the input sizes of the weights.  Therefore, the
    theorem follows by Lemmas \ref{hives} and \ref{bz}.
    
    For the second part of the theorem regarding type $A_r$, the hive
    polytopes provide a very fast method for determining whether \(
    c_{\lambda \mu}^{\nu} = 0 \).  According to the Saturation Theorem
    (see Section \ref{conjectures}), \( c_{\lambda \mu}^{\nu} = 0 \)
    if and only if the corresponding hive polytope is empty. Hence,
    it suffices to check whether the system of inequalities defining
    the hive polytope is feasible, which can be done in polynomial
    time for {\em arbitrary} dimension as a corollary of the
    polynomiality of linear programming via Khachian's ellipsoid
    algorithm (see \cite{schrijver}).

It is useful to notice that every Kostka number \( K_{\lambda \mu} \)
is a Littlewood--Richardson coefficient for some choice of highest
weights. One short bijection is given by \( K_{\lambda \mu} =
c_{\sigma \lambda}^{\tau} \), where
\[ 
    \begin{cases}
        \tau_{i}   = \mu_{i} + \mu_{i+1} + \dotsb \\
        \sigma_{i} = \mu_{i+1} + \mu_{i+2} + \dotsb
    \end{cases}
    \qquad
    \text{for \( i = 1, 2, \dotsc, r \)}.
\]

\begin{figure}[tbp]
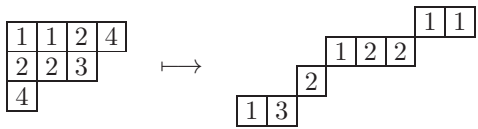

    \[ 
        \young(1124,223,4)
        \quad
        \longmapsto
        \quad
        \young(::::::11,:::122,::2,13)
    \]
     \caption{Corresponding semi-standard Young and
     Littlewood--Richardson tableaux}
    \label{fig:1}
\end{figure}

For those familiar with the enumeration of semi-standard Young
tableaux and Little\-wood--Rich\-ard\-son tableaux by Kostka numbers
and Little\-wood--Rich\-ard\-son coefficients respectively (see, {\em
e.g.}, \cite{stanleyII}), the bijection is straightforward: Given a
semi-standard Young tableau \( Y \) with shape \( \lambda \) and
content \( \mu \), construct a Littlewood--Richardson tableau \( L \)
with shape \( \tau/\sigma \) and content \( \lambda \) by filling the
boxes as follows.  Start with a skew Young diagram \( D \) with shape
\( \tau/\sigma \).  For \( j = 1, \dotsc, r \), and for each entry \(
i \) in the \( j \)th row of \( Y \), place a \( j \) in the \( i \)th
row of \( D \).  Let \( L \) be the tableau produced by filling the
boxes of \( D \) in this fashion.  (See Figure \ref{fig:1}.)  It is
not hard to see that, under this map, the column-strictness condition
on \( Y \) is equivalent to the lattice permutation condition on \( L
\).  It follows that the map just described is a bijection between
semi-standard Young tableaux with shape \( \lambda \) and content \(
\mu \) and Littlewood--Richardson tableaux with shape \( \tau/\sigma
\) and content \( \lambda \).  Thus, computing Kostka numbers reduces
to computing Littlewood--Richardson coefficients. 
 
As a corollary, we get a similar result for the computation of Kostka 
numbers based on polyhedral methods.

\begin{cor}
    For fixed rank $r$, the Kostka number $K_{\lambda \mu}$ can be
    computed in polynomial time in the size of the input weights
    $\lambda$ and $\mu$.  For arbitrary rank one can decide in
    polynomial time in the size of the weights whether $K_{\lambda
    \mu} = 0$ or not.
\end{cor}

\section{Using the Algorithm in Practice} \label{comparison}

Using the explicit definitions for the hive and BZ-polytopes as the
sets of solutions to systems of linear inequalities and equalities, we
wrote a \texttt{Maple} notebook which, when given a triple of highest
weights, produces the corresponding hive or BZ--polytope in
\texttt{LattE} readable input format.  The notebook is available from
\url{http://math.ucdavis.edu/~tmcal}.  All computations were done on a
Linux PC with a 2 GHZ CPU and 4 Gigabytes of memory.

From our experiments, we conclude that (1) The polyhedral method of
computing tensor product multiplicities complements the method
employed in \texttt{LiE}.  \texttt{LiE} is effective for larger ranks
(up to $r=10$, say), but the sizes of the weights must be kept small.
This is because \texttt{LiE} uses the Klimyk formula to generate the
entire direct sum decomposition of the tensor product, after which it
dispenses with all but the single desired term.  However, computing
all of the terms in the direct sum decomposition is not feasible when
the sizes of the entries in the weights grow into the 100s.  On the
other hand, (2) Lattice point enumeration is often effective for very
large weights (in particular, the algorithm is suitable for
investigating the stretching properties of Section \ref{conjectures}).
However, the rank must be relatively low (roughly $r \leq 6$) because
lattice point enumeration complexity grows exponentially in the
dimension of the polytope, and the dimensions of these polytopes grow
quadratically with the rank of the Lie algebra.  Together, the two
algorithms cover a larger range of problems.

We would also like to mention that other authors are also using
lattice points in polytopes to compute Clebsch--Gordan coefficients,
this time via the calculation of vector partition functions for
classical root systems \cite{bbcv, cochet1, cochet2}.  These authors
report that, like us, they can compute with large sizes of weight
entries.

\subsection{Experiments for type \( A_{r} \)}
\label{hivespolytopes}

In the tables below, we express highest weights in terms of the
canonical basis \( e_{1}, \dotsc, e_{r} \), so that the highest
weights are partitions with \( r \) parts.

Experiments indicate that lattice point enumeration is very efficient
for computing Littlewood--Richardson numbers when \( r \leq 5 \).
First, we computed over 30 instances with randomly generated weights
with leading entries larger than 40 with our approach and with
\texttt{LiE}.  In all cases our algorithm was faster.  After that, we
did a ``worst case'' sampling for the table in Figure
\ref{fig:Artable} comparing the computation times of \texttt{LattE}
and \texttt{LiE}.  To produce the \( i \)th row of that table, we
selected uniformly at random 1000 triples of weights \( (\lambda, \mu,
\nu) \) in which the largest parts of \( \lambda \) and \( \mu \) were
bounded above by \( 10i \) and \( \abs{\nu} = \abs{\lambda} +
\abs{\mu} \) (this is a necessary condition for \( c_{\lambda,
\mu}^{\nu} \ne 0 \)).  Then we evaluated the corresponding hive
polytopes with \texttt{LattE}.  The \texttt{LattE} input files are
created with our \texttt{Maple} program.  The weight triple in the \(
i \)th row is the one that \texttt{LattE} took the longest time to
compute.  We then computed the same tensor product multiplicity with
\texttt{LiE}.  The table in Figure \ref{fig:Artablelarge}, shows the
running time needed when using \texttt{LattE} to compute weight
triples with entries in the thousands or millions.

When \( r \geq 6 \), the running time under \texttt{LattE} begins to
blow up.  Still, for \( r = 6 \), all examples we attempted could be
computed in under 30 minutes using \texttt{LattE}, and most could be
computed in under 5 minutes.  For example, among 54 nonempty hive
polytopes chosen uniformly at random among those in which the weights
had entries less than 100, all but seven could be computed in under 5
minutes with \texttt{LattE}, and the remaining seven could all be
computed in under 30 minutes.  None of these Littlewood--Richardson
coefficients could be computed with \texttt{LiE}.  At \( r = 7 \),
lattice point enumeration becomes less effective, with examples
typically taking several hours or more to evaluate.

\begin{sidewaysfigure}
    \centering
    \begin{tabular}{|l|r|r|r|}
        \hline
        \( \lambda, \mu, \nu \)                              & \( c_{\lambda \mu}^{\nu} \) & \texttt{LattE} runtime & \texttt{LiE} runtime \\
        \hline
        \hline
        (9,7,3,0,0), (9,9,3,2,0), (10,9,9,8,6)               & 2                            & 0m00.74s                &  0m00.01s            \\
        (18,11,9,4,2), (20,17,9,4,0), (26,25,19,16,8)        & 453                          & 0m03.86s                &  0m00.12s            \\
        (30,24,17,10,2), (27,23,13,8,2), (47,36,33,29,11)    & 5231                         & 0m05.21s                &  0m02.71s            \\
        (38,27,14,4,2), (35,26,16,11,2), (58,49,29,26,13)    & 16784                        & 0m06.33s                &  0m25.31s            \\
        (47,44,25,12,10),(40,34,25,15,8),(77,68,55,31,29)    & 5449                         & 0m04.35s                &  1m55.83s            \\
        (60,35,19,12,10),(60,54,27,25,3),(96,83,61,42,23)    & 13637                        & 0m04.32s                & 23m32.10s            \\
        (64,30,27,17,9), (55,48,32,12,4), (84,75,66,49,24)   & 49307                        & 0m04.63s                & 45m52.61s            \\
        (73,58,41,21,4), (77,61,46,27,1), (124,117,71,52,45) & 557744                       & 0m07.02s                & \( > \) 24 hours     \\
        \hline
    \end{tabular}
    \caption{A sample comparison of running times between \texttt{LattE} and \texttt{LiE} case of $A_r$}
    \label{fig:Artable}
    
    \bigskip
    \bigskip
    
    \begin{tabular}{|l|r|r|}
        \hline
        \( \lambda, \mu, \nu \)                                            & \( c_{\lambda \mu}^{\nu} \) & \texttt{LattE} runtime \\
        \hline
        \hline
        \begin{tabular}{l} (935,639,283,75,48) \\ (921,683,386,136,21) \\ (1529,1142,743,488,225) \end{tabular} & 1303088213330                & 0m07.84s                \\
        \hline
        \begin{tabular}{l} (6797,5843,4136,2770,707) \\ (6071,5175,4035,1169,135) \\ (10527,9398,8040,5803,3070) \end{tabular} & 459072901240524338 & 0m09.63s         \\
        \hline
        \begin{tabular}{l} (859647,444276,283294,33686,24714) \\ (482907,437967,280801,79229,26997) \\ (1120207,699019,624861,351784,157647) \end{tabular} & 11711220003870071391294871475 & 0m08.15s \\
        \hline
    \end{tabular}
    \caption{Computing large weights with \texttt{LattE}. Case of $A_r$}
    \label{fig:Artablelarge}
\end{sidewaysfigure}

\subsection{Experiments for types \( B_{r} \), \( C_{r} \), and \( D_{r} \)}
\label{BZ-polytopes}

To compute Clebsch--Gordan coefficients in types \( B_{r} \), \( C_{r}
\), and \( D_{r} \), we used the BZ-polytopes.  In the tables that
follow, all weights are given in the basis of fundamental weights for
the corresponding Lie algebra.

Our experiments followed the same process we used for $A_r$: First,
for each root system, we computed over 30 instances with randomly
generated weights with entries larger than 40 with our approach and
with \texttt{LiE}.  In all cases our algorithm was faster.  After
that, we did a ``worst case'' sampling to produce the tables in Figure
\ref{fig:BCDrtable} comparing the computation times of \texttt{LattE}
and \texttt{LiE}.  As in Section \ref{hivespolytopes}, these weight
triples were the ones which \texttt{LattE} took the longest to
evaluate among thousands of instances generated with the following
procedure: First, to produce line $i$ of a table, we set an upper
bound $U_i$ for the entries of each weight.  Then, we generated 1000
random weight triples with entry sizes no larger than $U_i$.  Here are
the specific values of \( U_{i} \) used in each of the three tables in
Figure \ref{fig:BCDrtable}.  For type $B_r$, the bounds $U_i$ were 50,
60, 70, and 10,000, respectively.  For type $C_r$, the bounds $U_i$
were 50, 60, 80, and 10,000, respectively.  Finally, for type $D_r$,
the bounds $U_i$ were 20, 30, 40, and 10,000, respectively.  For each
generated triple of weights, we produced the associated BZ-polytopes
(using our \texttt{Maple} notebook) and counted their lattice points
with \texttt{LattE}.  In the table are those instances that were
slowest in \texttt{LattE}.  We also recorded in the table the time
taken by \texttt{LiE} for the same instances.  One can see the running
time needed by \texttt{LattE} is hardly affected by growth in the size
of the input weights, while the time needed by \texttt{LiE} grows
rapidly.

We found that for types \( B_{r} \) and \( C_{r} \), lattice point
enumeration with the BZ-polytopes is very effective when \( r \leq 3
\).  Each of the many thousands of examples we generated could be
evaluated by \texttt{LattE} in under 10 seconds (the examples in
Figure \ref{fig:BCDrtable} were the worst cases).  When \( r=4 \), the
computation time begins to blow up, with examples typically taking
half an hour or more to compute.  The polyhedral method is also
reasonably efficient for type \( D \) Lie algebras with rank 4, the
lowest rank in which they are defined.  All of the examples we
generated could be evaluated by \texttt{LattE} in under 5 minutes.

\begin{sidewaysfigure}
    \centering
    \begin{tabular}{|c|l|r|r|r|}
        \hline
                                     & \( \lambda, \mu, \nu \)            & \( C_{\lambda \nu}^{\mu} \) & \texttt{LattE} runtime & \texttt{LiE} runtime \\
        \hline
        \hline
        \multirow{4}{*}{\( B_{3} \)} & (46,42,38), (38,36,42), (41,36,44) & 354440672                   & 0m09.58s                & 1m45.27s             \\
                                     & (46,42,41), (14,58,17), (50,54,38) & 88429965                    & 0m06.38s                & 3m16.01s             \\
                                     & (15,60,67), (58,70,52), (57,38,63) & 626863031                   & 0m07.14s                & 6m01.43s             \\
                                     & (5567,2146,6241), (6932,1819,8227), (3538,4733,3648) & 215676881876569849679 & 0m7.07s & n/a \\
        \hline
        \multirow{4}{*}{\( C_{3} \)} & (25,42,22), (36,38,50), (31,33,48) & 87348857  & 0m07.48s & 0m17.21s \\
                                     & (34,56,36), (44,51,49), (37,51,54) & 606746767 & 0m08.42s & 2m57.27s \\
                                     & (39,64,58), (65,15,72), (70,41,44) & 519379044 & 0m07.63s & 8m00.35s \\
                                     & (5046,5267,7266), (7091,3228,9528), (9655,7698,2728) & 1578943284716032240384 & 0m07.66s & n/a \\
        \hline
        \multirow{4}{*}{\( D_{4} \)} & (13,20,10,14), (10,20,13,20), (5,11,15,18) & 41336415    & 2m46.88s & 0m12.29s \\
                                     & (12,22,9,30),(28,14,15,26),(10,24,10,26)   & 322610723   & 3m04.31s & 7m03.44s \\
                                     & (37,16,31,29),(40,18,35,41),(36,27,19,37)  & 18538329184 & 4m29.63s & \( > \)60m \\
                                     & (2883,8198,3874,5423),(1901,9609,889,4288),(5284,9031,2959,5527) & 1891293256704574356565149344 & 2m06.42s & n/a \\
        \hline
    \end{tabular}
    \caption{A sample comparison of running times between \texttt{LattE} and \texttt{LiE}}
    \label{fig:BCDrtable}
\end{sidewaysfigure}

\section{Two Conjectures that Could Generalize the Saturation Theorem}
\label{conjectures}

In 1998, Knutson and Tao used the hive polytopes to prove the
\emph{Saturation Theorem}.  Buch has written a very clear exposition
of this proof in \cite{buch}.

\begin{thm}
    [\cite{knutsontao}] (Saturation) Given highest weights \( \lambda
    \), \( \mu \), and \( \nu \) for a Lie algebra of type $A_r$, and
    given an integer $n > 0$, the Littlewood--Richardson coefficient
    $c_{\lambda \mu}^\nu$ satisfies
    \[ 
        c_{\lambda \mu}^{\nu} \neq 0 
        \quad \Longleftrightarrow \quad
        c_{n\lambda, n\mu}^{n\nu} \neq 0.
    \]
\end{thm}

In hive polytope language, the Saturation Theorem can be restated as
\[
    \# \left( H_{\lambda \mu}^{\nu} \cap \Z^{d} \right) \neq \varnothing
    \quad \Longleftrightarrow \quad
    \# \left( H_{n\lambda, n\mu}^{n\nu} \cap \Z^{d} \right) \neq \varnothing,
\]
where \( d = {r+2 \choose 2} \).  The definition of hive polytopes
(see Definition \ref{def:HivePolytopes} above) implies that \(
H_{n\lambda, n\mu}^{n\nu} = n\,H_{\lambda \mu}^{\nu} \), so the
Saturation Theorem is equivalent to the existence of a lattice point
inside each hive polytope \( H_{\lambda \mu}^{\nu} \). Now we would
like to state two conjectures that generalize this theorem.

\subsection{First Conjecture}

To show that every hive polytope contains an integral point, Knutson
and Tao actually proved that every hive polytope contains an integral
{\em vertex}.  Our idea was to take a different approach to show a
generalization of this last result using the basic theory of
triangulations of semigroups.  To develop this idea, observe that the
boundary equalities and rhombus inequalities that define a hive
polytope may be expressed as the set of solutions to a system of
matrix equalities and inequalities:
\begin{equation}\label{eq:HiveInequalities}
    H_{\lambda \mu}^{\nu}
    =
    \left\lbrace h \in \R^{(r+1)(r+2)/2} :
        \begin{array}{l}
            B h =    b(\lambda, \mu, \nu), \\
            R h \leq 0
        \end{array}
    \right\rbrace,
\end{equation}
where $B$ and $R$ are integral matrices (depending on \( r \)), and
$b(\lambda, \mu, \nu)$ is a integral vector depending on \( \lambda
\), \( \mu \), and \( \nu \).  Here we think of a hive pattern \( h \)
as a column vector of dimension \( (r+1)(r+2)/2 \).  Note that there
is some degree of choice in how the boundary equalities and rhombus
inequalities are encoded as matrices \( B \) and \( R \),
respectively.  However, all such encodings are equivalent for our
purposes.

A polytope defined by such a system of equalities and inequalities may
be {\em homogenized} by adding ``slack variables''.  This produces an
equivalent polytope defined as the set of nonnegative solutions to a
system linear equalities.  Following this procedure, we define the
{\em homogenized hive polytope} \( \tilde{H}_{\lambda \mu}^{\nu} \) by
\[
    \tilde{H}_{\lambda, \mu}^{\nu}
    =
    \left\lbrace
        \tilde{h}
        :
        \left[
            \begin{matrix}
                B & 0 \\
                R & I
            \end{matrix}
        \right]
        \tilde{h}
        =
        \left[
            \begin{matrix}
                b(\lambda, \mu, \nu) \\
                0
            \end{matrix}
        \right],
        \quad
        \tilde{h} \geq 0
    \right\rbrace
\]
(where \( I \) is the identity matrix).  The equivalence between \(
H_{\lambda \mu}^{\nu} \) and \( \tilde{H}_{\lambda \mu}^{\nu} \) is
given by the linear map
\[
    h
    \mapsto
    \left[
        \begin{matrix}
            h    \\
            -R h
        \end{matrix}
    \right].
\]
Note that this linear map preserves vertices and integrality.
Therefore, to prove the Saturation Theorem, it suffices to show that
every homogenized hive polytope contains an integral vertex. 
Proceeding with this idea, we make the following definitions.

\begin{definition}
    Fix \( r \in \Z \).  Define the \emph{homogenized hive matrix} to
    be
    \[
        M
        =
        \left[
            \begin{matrix}
                B & 0 \\
                R & I
            \end{matrix}
        \right]
    \]
    (where \( B \) and \( R \) are as in equation
    \eqref{eq:HiveInequalities}).  Given an integral vector \( b \)
    with dimension equal to the number of rows in \( M \), define the
    {\em generalized hive polytope} or {\em g-hive polytope} \( H_{b}
    \) by
    \begin{equation} \label{ghives}
        H_b
        =
        \left\lbrace
            \tilde{h}
            :
            M
            \tilde{h}
            =
            b,
            \,
            \tilde{h} \geq 0
        \right\rbrace.
    \end{equation}
\end{definition}

Note that the homogenized hive polytopes are g-hive polytopes which
satisfy additional restrictions to the right-hand side vector (such
has the final entries of \( b \) being 0).

We now state some very basic facts about vertices of polyhedra
expressed in the form \( \{x : A x = b,\, x \geq 0\} \).  Let a finite
collection of integral points \( \{ a_{1}, \dotsc, a_{n} \} \subset
\Z^{m} \) be given, and let \( A \) be the matrix with columns \(
a_{1}, \dotsc, a_{n} \).  Define \( \cone A \) to be the cone in \(
\R^{m} \) generated by the point-set \( \{ a_{1}, \dotsc, a_{n} \} \):
\[
    \cone A
    =
    \{x_{1}a_{1} + \dotsb + x_{n}a_{n} \, : \, x_{1}, \dotsc, x_{n} \geq 0 \}.
\]
Then, for each vector \( b\in \Z^{m} \), we have a polytope
\[
    P_{b} = \{x : Ax = b,\, x \geq 0 \},
\]
and \( P_{b} \neq \varnothing \) if and only if \( b \in \cone A \).
In other words, there is a correspondence between nonempty polytopes
\( P_{b} \), \( b \in \Z^{m} \), and the elements of the semigroup
generated by the columns of \( A \).  The crucial property for our
purposes is the following.

\begin{lemma}
    If \( b \in (\cone A') \cap \Z^{m} \) for some \( m \times m \)
    submatrix \( A' \) of \( A \) with \( \det A' = \pm 1 \), then \(
    P_{b} \) has an integral vertex.
\end{lemma}

\begin{proof}
    Suppose that \( b \in (\cone A') \cap \Z^{m} \) for some \( m
    \times m \) submatrix \( A' \) of \( A \) with \( \det A' = \pm 1
    \).
    
    Let the columns of \( A' \) be \( a_{i_{1}}, \dotsc, a_{i_{m}} \),
    and let \( J = \{i_{1}, \dotsc, i_{m}\} \) be the indices of these
    columns.  Then there is a vector \( x = (x_{1}, \dotsc, x_{n})^{T}
    \in \R^{n}_{\geq 0} \) such that \( Ax = b \) and \( x_{i} = 0 \)
    for each \( i \notin J \).  Letting \( x' = (x_{i_{1}}, \dotsc,
    x_{i_{m}}) \) and using Crammer's rule to solve for \( x' \) in \(
    A'x' = b \), we find that \( x \) is an integral vector.  Thus, \(
    x \) is an integral lattice point in the polytope \( P_{b} \).
    
    To see that \( x \) is in fact a vertex of \( P_{b} \), Recall
    that the codimension (with respect to the ambient space) of the
    face containing a solution to a system of linear equalities and
    inequalities is the number of linearly independent equalities or
    inequalities satisfied with {\em equality}.  Observe that \( x \)
    is a solution to the system of \( n \) equalities
    \[ 
        \begin{cases}
            Ax = b     &             \\
            x_{i} = 0, & i \notin J.
        \end{cases}
    \]
    We claim that this is a linearly independent system of equalities.
    For suppose otherwise.  Then the zero vector is a nontrivial
    linear combination of the rows of \( A \) and the row-vectors \(
    e_{i} \), \( i \notin J \).  But this implies that the zero vector
    is a nontrivial linear combination of the rows of \( A' \), which
    is impossible because \( \det A' \ne 0 \).
    
    Thus, having shown that \( x \) satisfies the \( n \) linearly
    independent equalities above, we have shown that \( x \) lies in a
    codimension-\( n \) face of \( P_{b} \), {\em i.e.}, \( x \) is a
    vertex.
\end{proof}

We say that \( a_{i_{1}}, \dotsc, a_{i_{m}} \) is a {\em unimodular}
subset if the submatrix \( A' \) of \( A \) with columns \( a_{i_{1}},
\dotsc, a_{i_{m}} \) satisfies \( \det A' = \pm 1 \).  We say that the
matrix \( A \) has a {\em unimodular cover} (resp.  {\em unimodular
triangulation}) if the point set \( \{ a_{1}, \dotsc, a_{n} \} \) has
a unimodular cover (resp.  unimodular triangulation).

\begin{cor}
    If \( A \) has a unimodular cover, then \( P_{b} \) has an
    integral vertex for every integral \( b \in \cone(A)\).
\end{cor}

Our conjecture is that this corollary applies to the homogenized hive
matrix.  More precisely, we conjecture the following.

\begin{conjecture}{\label{conj:conj1}}        
    For each \( r \), the homogenized hive matrix has a unimodular
    triangulation.  Consequently, every g-hive polytope has an
    integral vertex.
\end{conjecture}

Since the hive polytopes are special cases of the g-hive polytopes,
Conjecture \ref{conj:conj1} generalizes the Saturation Theorem.

\begin{thm}\label{unimodthm}
    Conjecture (\ref{conj:conj1}) is true for \( r \leq 6 \).
\end{thm}

To compute the unimodular triangulations that provide a proof of
Theorem \ref{unimodthm} we used the software {\tt topcom}
\cite{topcom}.  It may be worth noting that the triangulations used to
prove Theorem \ref{unimodthm} were all placing triangulations.

\subsection{Second Conjecture}

For our second conjecture, we looked at general properties satisfied
by the tensor product multiplicities for semisimple Lie Algebras of
types \( B_r \), \( C_r \), and \( D_r \) under the operation of {\em
stretching of multiplicities}.  By {\em stretching of multiplicities},
we refer to the function $e \colon \Z_{>0} \rightarrow Z_{\geq 0}$
defined by \( e(n) = C^{n\nu}_{n\lambda, n\mu} \).

It follows from the definitions of the BZ-polytopes that, given any
highest weights \( \lambda, \mu, \nu \) of a semisimple Lie algebra,
\( C_{n\lambda, n\mu}^{n\nu} = e(n) \) is a quasi-polynomial in \( n
\).  Indeed, \( e(n) \) is, in polyhedral language, the \emph{Ehrhart
quasi-polynomial} of the corresponding BZ-polytope.  We recall the
basic theory of Ehrhart quasi-polynomials.  Its origins can be traced
to the work of Ehrhart \cite{ehrhart} in the 1960's (see Chapter 4 of
\cite{stanleyI} for an excellent introduction).

Given a convex polytope $P$, let
\[
    n P = \{x : (1/n)x \in P \}, \qquad n = 1,2,\dotsc. 
\]
If $P$ is a \( d \)-dimensional rational polytope in \( \R^{k} \),
then the counting function \( i_P (n) = \# (n P \cap \Z^k) \) is a
quasi-polynomial function of degree $d$; that is, there are
polynomials $f_1(n),\dotsc, f_N(n)$ of degree $d$ s.t.
\[
    i_P(n) 
    =
    \begin{cases}
        f_1(n) & \text{if $n \equiv 1 \quad \text{mod} \, N$}, \\
        \vdots &                                               \\
        f_N(n) & \text{if $n \equiv N \quad \text{mod} \, N$}.
    \end{cases}
\]

If we put \( P = H_{\lambda \mu}^{\nu} \), then the Ehrhart
quasi-poly\-nomial of \( P \) is just the stretched
Little\-wood--Rich\-ard\-son coefficient \( c_{n\lambda, n\mu}^{n\nu}
\).  The Ehrhart quasi-poly\-no\-mials of hive polytopes have been
studied by several authors.  Our experiments support the conjecture of
King et al \cite{kingetal} that the coefficients are in fact positive.
Since large weights can be computed with lattice point enumeration, it
is possible to produce the Ehrhart quasi-polynomials for the stretched
Clebsch--Gordan coefficients in the other types.  See Figures
\ref{fig:stretchedB3}--\ref{fig:stretchedD3} for some sample examples
out of the many hundreds generated.  Our experiments motivate the
following ``stretching conjecture''.

\begin{conjecture} \label{conjtwo}
    (Stretching Conjecture) Given highest weights \( \lambda,\mu, \nu \)
    of a Lie algebra of type \(A_{r}\), \( B_{r} \), \( C_{r} \), or \(
    D_{r} \), let
    \[
        C_{n\lambda, n\mu}^{n\nu}
        =
        \begin{cases}
            f_1(n) & \text{if $n \equiv 1 \quad \text{mod} \, N$}, \\
            \vdots &                                               \\
            f_N(n) & \text{if $n \equiv N \quad \text{mod} \, N$}
        \end{cases}
    \]
    be the quasi-polynomial representation of the stretched
    Clebsch--Gordan coefficient $C_{n\lambda, n\mu}^{n\nu}$.  Then the
    coefficients of each polynomial $f_{i}$ are all nonnegative.
\end{conjecture}

The type $A_r$ case of this conjecture was made by King, Tollu, and
Toumazet in \cite{kingetal}.  That Conjecture \ref{conjtwo} implies the
Saturation Theorem follows from a result of Derksen and Weyman
\cite{derksenweyman} showing that the Ehrhart quasi-polynomials of
Hive polytopes are in fact just polynomials.

We should remark that the saturation property is known not to hold in
the root systems $B_r$, $C_r$, and $D_r$.  A simple example in $B_2$,
due to Kapovich, Leeb, and Millson \cite{kapovichetal}, is given by
setting $\lambda = \mu = \nu = (1, 0)$ (with respect to the basis of
fundamental weights).  In this case we have
\[
    C_{n\lambda, n\mu}^{n\nu}
    =
    \begin{cases}
        0      & \text{if $n$ is odd}, \\
        1      & \text{if $n$ is even}.
    \end{cases}
\]
This example also demonstrates why the Stretching Conjecture is not
contradicted by the failure of the saturation property in the root
systems $B_r$, $C_r$, or $D_r$.  Since the stretched multiplicities
are not necessarily polynomials in these cases, it is possible for
them to evaluate to zero for some nonnegative integer while still
having all nonnegative coefficients.

\noindent {\bf Acknowledgements:} We are thankful to Arkady
Berenstein, Charles Cochet, Misha Kapovich, Allen Knutson, Peter
Littelmann, Ezra Miller, Hariharan Narayanan, J\"org Rambau, Francisco
Santos, Anna Schilling, Etienne Rassart, Mich\`ele Vergne, and Andrei
Zelevinsky for useful comments.

\newpage

\begin{sidewaysfigure}[tbp]
    \begin{tabular}{|l|l|r|}
        \hline
        \( \lambda, \mu, \nu \) & $C_{n \lambda, n \nu}^{n \mu}$ & \texttt{LattE} runtime \\
        \hline
        \hline
        \begin{tabular}{l}(0, 15, 5) \\ (12, 15, 3) \\ (6, 15, 6) \end{tabular} & \(\begin{cases} {\frac {68339}{64}}\,{n}^{5}+{\frac {407513}{384}}\,{n}^{4}+{\frac {13405}{32}}\,{n}^{3}+{\frac {9499}{96}}\,{n}^{2}+{\frac {107}{8}}\,n+1, \text{ \( n \) even} \\ {\frac {68339}{64}}\,{n}^{5}+{\frac {407513}{384}}\,{n}^{4}+{\frac {13405}{32}}\,{n}^{3}+{\frac {16355}{192}}\,{n}^{2}+{\frac {659}{64}}\,n+{\frac {75}{128}}, \text{ \( n \) odd}\end{cases}\) & 7m08.57s \\
        \hline
        \begin{tabular}{l}(4, 8, 11) \\ (3, 15, 10) \\ (10, 1, 3) \end{tabular} & \(\begin{cases} {\frac {13}{4}}\,{n}^{2}+3\,n+1, \text{ \( n \) even} \\ {\frac {13}{4}}\,{n}^{2}+3\,n+3/4, \text{ \( n \) odd}\end{cases}\) & 0m01.64s \\
        \hline
        \begin{tabular}{l}(8, 1, 3) \\ (11, 13, 3) \\ (8, 6, 14) \end{tabular} & \(\begin{cases} {\frac {121}{576}}\,{n}^{6}+{\frac {1129}{640}}\,{n}^{5}+{\frac {6809}{1152}}\,{n}^{4}+{\frac {163}{16}}\,{n}^{3}+{\frac {2771}{288}}\,{n}^{2}+{\frac {191}{40}}\,n+1, \text{ \( n \) even} \\ {\frac {121}{576}}\,{n}^{6}+{\frac {1129}{640}}\,{n}^{5}+{\frac {6809}{1152}}\,{n}^{4}+{\frac {1933}{192}}\,{n}^{3}+{\frac {659}{72}}\,{n}^{2}+{\frac {8003}{1920}}\,n+{\frac {93}{128}}, \text{ \( n \) odd}\end{cases}\) & 0m10.26s \\
        \hline
        \begin{tabular}{l}(8, 9, 14) \\ (8, 4, 5) \\ (1, 5, 15) \end{tabular} & \(\begin{cases} {\frac {4117}{192}}\,{n}^{6}+{\frac {50369}{640}}\,{n}^{5}+{\frac {14829}{128}}\,{n}^{4}+{\frac {703}{8}}\,{n}^{3}+{\frac {3541}{96}}\,{n}^{2}+{\frac {341}{40}}\,n+1, \text{ \( n \) even} \\ {\frac {4117}{192}}\,{n}^{6}+{\frac {50369}{640}}\,{n}^{5}+{\frac {14829}{128}}\,{n}^{4}+{\frac {5599}{64}}\,{n}^{3}+{\frac {3451}{96}}\,{n}^{2}+{\frac {5001}{640}}\,n+{\frac {97}{128}}, \text{ \( n \) odd}\end{cases}\) & 0m13.29s \\
        \hline
        \begin{tabular}{l}(10, 5, 6) \\ (5, 4, 10) \\ (0, 7, 12) \end{tabular} & \(\begin{cases} {\frac {669989}{960}}\,{n}^{5}+{\frac {286355}{384}}\,{n}^{4}+{\frac {10803}{32}}\,{n}^{3}+{\frac {7993}{96}}\,{n}^{2}+{\frac {1427}{120}}\,n+1, \text{ \( n \) even} \\ {\frac {669989}{960}}\,{n}^{5}+{\frac {286355}{384}}\,{n}^{4}+{\frac {10803}{32}}\,{n}^{3}+{\frac {15509}{192}}\,{n}^{2}+{\frac {10081}{960}}\,n+{\frac {65}{128}}, \text{ \( n \) odd}\end{cases}\) & 2m52.39s \\
        \hline
    \end{tabular}
    \caption{Stretched Clebsch--Gordan coefficients for \( B_{3} \).}
    \label{fig:stretchedB3}
\end{sidewaysfigure}

\begin{sidewaysfigure}[tbp]
    \begin{tabular}{|l|l|r|}
        \hline
        \( \lambda, \mu, \nu \) & $C_{n \lambda, n \nu}^{n \mu}$ & \texttt{LattE} runtime \\
        \hline
        \hline
        \begin{tabular}{l} (1,13,6) \\ (14,15,5) \\ (5,11,7) \end{tabular} & \(\begin{cases} {\frac {5937739}{5760}}\,{n}^{6}+{\frac {87023}{40}}\,{n}^{5}+{\frac {936097}{576}}\,{n}^{4}+{\frac {27961}{48}}\,{n}^{3}+{\frac {85397}{720}}\,{n}^{2}+{\frac {883}{60}}\,n+1 , \text{ \( n \) even} \\ {\frac {5937739}{5760}}\,{n}^{6}+{\frac {87023}{40}}\,{n}^{5}+{\frac {936097}{576}}\,{n}^{4}+{\frac {27961}{48}}\,{n}^{3}+{\frac {657931}{5760}}\,{n}^{2}+{\frac {3097}{240}}\,n+3/4, \text{ \( n \) odd}\end{cases}\) & 21m20.59s \\
        \hline
        \begin{tabular}{l} (4,15,14) \\ (12,12,10) \\ (4,9,8) \end{tabular} & \(\begin{cases} {\frac {22199219}{2880}}\,{n}^{6}+{\frac {8154617}{960}}\,{n}^{5}+{\frac {4500665}{1152}}\,{n}^{4}+{\frac {31297}{32}}\,{n}^{3}+{\frac {226903}{1440}}\,{n}^{2}+{\frac {2021}{120}}\,n+1, \text{ \( n \) even} \\ {\frac {22199219}{2880}}\,{n}^{6}+{\frac {8154617}{960}}\,{n}^{5}+{\frac {4500665}{1152}}\,{n}^{4}+{\frac {31297}{32}}\,{n}^{3}+{\frac {217363}{1440}}\,{n}^{2}+{\frac {13513}{960}}\,n+{\frac {85}{128}}, \text{ \( n \) odd}\end{cases}\) & 17m05.74s \\
        \hline
        \begin{tabular}{l} (9,0,8) \\ (8,12,9) \\ (7,7,3) \end{tabular} & \(1/30\,{n}^{5}+3/8\,{n}^{4}+{\frac {19}{12}}\,{n}^{3}+{\frac {25}{8}}\,{n}^{2}+{\frac {173}{60}}\,n+1\) & 0m00.61s \\
        \hline
        \begin{tabular}{l} (10,2,7) \\ (8,10,1) \\ (7,5,5) \end{tabular} & \(\begin{cases} {\frac {596153}{1152}}\,{n}^{6}+{\frac {53425}{48}}\,{n}^{5}+{\frac {502621}{576}}\,{n}^{4}+{\frac {5577}{16}}\,{n}^{3}+{\frac {11941}{144}}\,{n}^{2}+{\frac {149}{12}}\,n+1, \text{ \( n \) even} \\ {\frac {596153}{1152}}\,{n}^{6}+{\frac {53425}{48}}\,{n}^{5}+{\frac {502621}{576}}\,{n}^{4}+{\frac {5577}{16}}\,{n}^{3}+{\frac {94097}{1152}}\,{n}^{2}+{\frac {131}{12}}\,n+{\frac {23}{32}}, \text{ \( n \) odd}\end{cases}\) & 19m24.55s \\
        \hline
        \begin{tabular}{l} (10,10,15) \\ (11,3,15) \\ (10,7,15) \end{tabular} & \(\begin{cases} {\frac {6084163}{320}}\,{n}^{6}+{\frac {507527}{30}}\,{n}^{5}+{\frac {1185853}{192}}\,{n}^{4}+{\frac {59995}{48}}\,{n}^{3}+{\frac {43039}{240}}\,{n}^{2}+{\frac {357}{20}}\,n+1, \text{ \( n \) even} \\ {\frac {6084163}{320}}\,{n}^{6}+{\frac {507527}{30}}\,{n}^{5}+{\frac {1185853}{192}}\,{n}^{4}+{\frac {59995}{48}}\,{n}^{3}+{\frac {144751}{960}}\,{n}^{2}+{\frac {883}{80}}\,n+{\frac {25}{64}}, \text{ \( n \) odd}\end{cases}\) & 16m05.08s \\
        \hline
    \end{tabular}
    \caption{Stretched Clebsch--Gordan coefficients for \( C_{3} \).}
    \label{fig:stretchedC3}
\end{sidewaysfigure}

\begin{sidewaysfigure}[tbp]
    \begin{tabular}{|l|l|r|}
        \hline
        \( \lambda, \mu, \nu \) & $C_{n \lambda, n \nu}^{n \mu}$ & \texttt{LattE} runtime \\
        \hline
        \hline
        \begin{tabular}{l}(0, 2, 10, 5) \\ (4, 11, 9, 11) \\ (5, 8, 6, 9) \end{tabular} & \(\begin{cases}{\frac {625007}{10080}}\,{n}^{7}+{\frac {729157}{2880}}\,{n}^{6}+{\frac {77197}{180}}\,{n}^{5}+{\frac {449539}{1152}}\,{n}^{4}+{\frac {298979}{1440}}\,{n}^{3}+{\frac {95189}{1440}}\,{n}^{2}+{\frac {10079}{840}}\,n+1, \text{ \( n \) even} \\ {\frac {625007}{10080}}\,{n}^{7}+{\frac {729157}{2880}}\,{n}^{6}+{\frac {77197}{180}}\,{n}^{5}+{\frac {449539}{1152}}\,{n}^{4}+{\frac {298979}{1440}}\,{n}^{3}+{\frac {95189}{1440}}\,{n}^{2}+{\frac {10079}{840}}\,n+{\frac {127}{128}}, \text{ \( n \) odd}\end{cases}\) & 20m24.79s  \\
        \hline
        \begin{tabular}{l}(2, 7, 12, 2) \\ (11, 10, 5, 9) \\ (13, 11, 1, 1) \end{tabular} & \(\begin{cases} {\frac {34675903}{80640}}\,{n}^{8}+{\frac {3037051}{1680}}\,{n}^{7}+{\frac {9121453}{2880}}\,{n}^{6}+{\frac {241181}{80}}\,{n}^{5}+{\frac {615083}{360}}\,{n}^{4}+{\frac {8947}{15}}\,{n}^{3}+{\frac {107791}{840}}\,{n}^{2}+{\frac {6721}{420}}\,n+1 \text{ \( n \) even} \\ {\frac {34675903}{80640}}\,{n}^{8}+{\frac {3037051}{1680}}\,{n}^{7}+{\frac {9121453}{2880}}\,{n}^{6}+{\frac {241181}{80}}\,{n}^{5}+{\frac {615083}{360}}\,{n}^{4}+{\frac {8947}{15}}\,{n}^{3}+{\frac {107791}{840}}\,{n}^{2}+{\frac {6721}{420}}\,n+{\frac {239}{256}} \text{ \( n \) odd}\end{cases}\) & 123m59.76s \\
        \hline
        \begin{tabular}{l}(3, 11, 0, 10) \\ (2, 15, 10, 15) \\ (10, 12, 11, 0) \end{tabular} & \(\begin{cases} {\frac {53609}{60}}\,{n}^{6}+{\frac {25631}{15}}\,{n}^{5}+{\frac {63779}{48}}\,{n}^{4}+{\frac {1627}{3}}\,{n}^{3}+{\frac {2497}{20}}\,{n}^{2}+{\frac {239}{15}}\,n+1, \text{ \( n \) even} \\ {\frac {53609}{60}}\,{n}^{6}+{\frac {25631}{15}}\,{n}^{5}+{\frac {63779}{48}}\,{n}^{4}+{\frac {1627}{3}}\,{n}^{3}+{\frac {2497}{20}}\,{n}^{2}+{\frac {239}{15}}\,n+{\frac {15}{16}}, \text{ \( n \) odd}\end{cases}\) & 2m37.73s \\
        \hline
        \begin{tabular}{l}(10, 1, 12, 4) \\ (1, 12, 0, 3) \\ (0, 5, 3, 4) \end{tabular} & \( 5\,{n}^{2}+4\,n+1 \) & 0m01.63s \\
        \hline
        \begin{tabular}{l}(12, 2, 5, 13) \\ (15, 6, 10, 11) \\ (2, 0, 12, 13) \end{tabular} & \(\begin{cases} {\frac {455263}{2016}}\,{n}^{7}+{\frac {447281}{576}}\,{n}^{6}+{\frac {198433}{180}}\,{n}^{5}+{\frac {971011}{1152}}\,{n}^{4}+{\frac {108787}{288}}\,{n}^{3}+{\frac {28969}{288}}\,{n}^{2}+{\frac {12631}{840}}\,n+1, \text{ \( n \) even} \\ {\frac {455263}{2016}}\,{n}^{7}+{\frac {447281}{576}}\,{n}^{6}+{\frac {198433}{180}}\,{n}^{5}+{\frac {971011}{1152}}\,{n}^{4}+{\frac {108787}{288}}\,{n}^{3}+{\frac {28969}{288}}\,{n}^{2}+{\frac {12631}{840}}\,n+{\frac {127}{128}}, \text{ \( n \) odd}\end{cases}\) & 4m25.90s \\
        \hline
    \end{tabular}
    \caption{Stretched Clebsch--Gordan coefficients for \( D_{4} \).}
    \label{fig:stretchedD3}
\end{sidewaysfigure}

\end{document}